\documentclass[12pt]{article}
\usepackage{amssymb}
\usepackage{amsmath}

\newtheorem{theorem}{Théorème}

\newtheorem{corollary}[theorem]{Corollaire}

\newtheorem{definition}[theorem]{Définition}

\newtheorem{lemma}[theorem]{Lemme}

\newtheorem{proposition}[theorem]{Proposition}
\newtheorem{remark}[theorem]{Remarque}

\input{tcilatex}

\begin{document}

\title{Quelques remarques sur les espaces d'interpolation $A^{\beta }$ }
\author{Daher Mohammad \\
%EndAName
D\'{e}partement de Maths, universit\'{e} Paris VII\\
e-mail: m.daher@orange.fr}
\maketitle

\noindent\ \ \ \ \ \ \emph{R\'{e}sum\'{e}}: Soient $(A_{0},A_{1})$ un couple
d'interpolation et $B_{j}$ l'adh\'{e}rence de $A_{0}^{\ast }\cap A_{1}^{\ast
}$ dans $A_{j}^{\ast },$ $j=0,1.$ Pour tout $\theta \in ]0,1[,$ il existe
une contraction injective naturelle $R^{\theta }:A^{\theta }\rightarrow
(B_{0}^{\ast },B_{1}^{\ast })^{\theta }.$ Dans ce travail on montre que $%
A^{\theta }=A_{\theta },$ pour tout $\theta \in \left] 0,1\right[ $ si $%
Z^{\beta }$ v\'{e}rifie quelques hypoth\`{e}ses raisonnables, pour un $\beta
\in ]0,1[,$ o\`{u} $Z^{\beta }$ est l'adh\'{e}rence de $R^{\beta }(A^{\beta
})$ dans $(B_{0}^{\ast },B_{1}^{\ast })^{\beta }.$

\emph{Abstract}: Let $(A_{0},A_{1})$ be an interpolation couple and let $%
B_{j}$ be the closure of $A_{0}^{\ast }\cap A_{1}^{\ast }$ in $A_{j}^{\ast
}, $ $j=0,1.$ For every $\theta \in ]0,1[,$ there exists a natural one to
one contraction $R^{\theta }:A^{\theta }\rightarrow (B_{0}^{\ast
},B_{1}^{\ast })^{\theta }.$ We show that $A^{\theta }=A_{\theta }$ for
every $\theta \in ]0,1[$, if $Z^{\beta }$ satisfies some reasonable
assumtions, for some $\beta \in ]0,1[,$ here $Z^{\beta }$ is the closure of $%
R^{\beta }(A^{\beta }) $ in $(B_{0}^{\ast },B_{1}^{\ast })^{\beta }$.


\emph{Mots cl\'{e}s}: Interpolation

\begin{center}
\medskip
\end{center}

\section{Introduction et notations}

On note $X^{\ast }$ le dual d'un espace de Banach $X.$

Soit $\overline{A}=(A_{0},A_{1})$ un couple d'interpolation complexe, au
sens de \cite{BL}. Soit $S=\left\{ z\in \mathbb{C};\text{ }0\leq \func{Re}%
(z)\leq 1\right\} .$

Rappelons d'abord la d\'{e}finition de l'espace d'interpolation $A_{\theta
}, $ o\`{u} $\theta \in \left] 0,1\right[ $ \cite[chap.4]{BL}. On note $%
\mathcal{F}\mathfrak{(}\overline{A})$ l'espace des fonctions $F$ \`{a}
valeurs dans $A_{0}+A_{1}$, continues born\'{e}es sur $S,$ holomorphes \`{a}
l'interieur de $S,$ telles que, pour $j\in \{0,1\},$ l'application $\tau \in 
\mathbb{R}\rightarrow F(j+i\tau )\in A_{j}$ est continue (prend ses valeurs
dans $A_{j})$ et $\left\Vert F(j+i\tau )\right\Vert _{A_{j}}\rightarrow 0,$
quand $\left\vert \tau \right\vert $ $\rightarrow +\infty .$ On munit $%
\mathcal{F}\mathfrak{(}\overline{A})$ de la norme

\begin{equation*}
\left\Vert F\right\Vert _{\mathcal{F}\mathfrak{(}\overline{A)}\text{ }}=\max
(sup_{\tau \in \mathbb{R}}\left\Vert F(i\tau )\right\Vert _{A_{0}},sup_{\tau
\in \mathbb{R}}\left\Vert F(1+i\tau )\right\Vert _{A_{1}}).
\end{equation*}

\noindent L'espace $A_{\theta }=(A_{0},A_{1})_{\theta }=\left\{ F(\theta );%
\text{ }F\in \mathcal{F}\mathfrak{(}\overline{A})\right\} $ est un Banach 
\cite[theorem 4.1.2]{BL} pour la norme d\'{e}finie par

\begin{equation*}
\left\Vert a\right\Vert _{A_{\theta }}=\inf \left\{ \left\Vert F\right\Vert
_{\mathcal{F}\mathfrak{(}\overline{A})};\text{ }F(\theta )=a\right\} .
\end{equation*}

Rappelons maintenant la d\'{e}finition de l'espace d'interpolation $%
A^{\theta }$ \cite[chapitre 4]{BL}. On note $\mathcal{G(}\overline{A})$
l'espace des fonctions $g$ \`{a} valeurs dans $A_{0}+A_{1}$, continues sur $%
S,$ holomorphes \`{a} l'interieur de $S,$ telles que $z\rightarrow
(1+\left\vert z\right\vert )^{-1}\left\Vert g(z)\right\Vert _{A_{0}+A_{1}}$
est born\'{e}e sur $S,$ $g(j+i\tau )-g(j+i\tau ^{\prime })\in A_{j}$ $\ $%
pour tous $\tau ,\tau ^{\prime }\in \mathbb{R}$, $j\in \{0,1\},$ et la
quantit\'{e} suivante est finie:%
\begin{eqnarray*}
&&\left\Vert g^{\cdot }\right\Vert _{Q\mathcal{G(}\overline{A}\mathcal{)}} \\
&=&\max \left[ \underset{\tau \neq \tau ^{\prime }}{\sup_{\tau ,\tau
^{\prime }\in R}}\left\Vert \frac{g(i\tau )-g(i\tau ^{\prime })}{\tau -\tau
^{\prime }}\right\Vert _{A_{0}},\underset{}{\sup_{\underset{\tau \neq \tau
^{\prime }}{\tau ,\tau ^{\prime }\in R}}}\left\Vert \frac{g(1+i\tau
)-g(1+i\tau ^{\prime })}{\tau -\tau ^{\prime }}\right\Vert _{A_{1}}\right] .
\end{eqnarray*}

\noindent Cette quantit\'{e} d\'{e}finit une norme sur l'espace $Q\mathcal{G(%
}\overline{A}\mathcal{)},$ quotient de $\mathcal{G(}\overline{A}\mathcal{)}$
par les applications constantes \`{a} valeurs dans $A_{0}\cap A_{1}$, et $Q%
\mathcal{G(}\overline{A}\mathcal{)}$ est complet pour cette norme \cite[%
lemma 4.1.3]{BL}.

On rappelle \cite[p. 89]{BL} que, pour $g\in \mathcal{G(}\overline{A}%
\mathcal{)},$

\begin{equation}
\left\Vert g^{\prime }(z)\right\Vert _{A_{0}+A_{1}}\leq \left\Vert g^{\cdot
}\right\Vert _{Q\mathcal{G}(\overline{A})},\;\;z\in S.  \label{2}
\end{equation}

\noindent C'est une cons\'{e}quence imm\'{e}diate de l'in\'{e}galit\'{e}

\begin{equation*}
\left\Vert \frac{g(z+it)-g(z)}{t}\right\Vert _{A_{0}+A_{1}}\leq \left\Vert
g^{\cdot }\right\Vert _{Q\mathcal{G}(\overline{A})},\;z\in S,\;t\in \mathbb{R%
}^{\ast },
\end{equation*}

\noindent qui d\'{e}coule de la d\'{e}finition de $\left\Vert g^{\cdot
}\right\Vert _{Q\mathcal{G}(\overline{A})}$ et du th\'{e}or\`{e}me des trois
droites \cite[lemma 1.1.2]{BL} appliqu\'{e} aux fonctions $z\rightarrow
\left\langle \frac{g(z+it)-g(z)}{t},a^{\ast }\right\rangle ,$ $t$ r\'{e}el
fix\'{e}, o\`{u} $a^{\ast }$ parcourt la boule unit\'{e} de $A_{0}^{\ast
}\cap A_{1}^{\ast }$.

L'espace $A^{\theta }=\left\{ \text{ }g^{\prime }(\theta );\text{ }g\in 
\mathcal{G(}\overline{A})\text{ }\right\} $ est un Banach \cite[theorem 4.1.4%
]{BL} pour la norme d\'{e}finie par

\begin{equation*}
\left\Vert a\right\Vert _{A^{\theta }}=\inf \left\{ \left\Vert g^{\cdot
}\right\Vert _{Q\mathcal{G}(\overline{A})};\text{ }g^{\prime }(\theta
)=a\right\} .
\end{equation*}

\noindent D'apr\`{e}s (\ref{2}) $\left\Vert a\right\Vert _{A_{0}+A_{1}}\leq
\left\Vert a\right\Vert _{A^{\theta }}.$ La contraction $A^{\theta
}\rightarrow A_{0}+A_{1}$ est injective par d\'{e}finition de $A^{\theta }.$

\noindent D'apr\`{e}s \cite{B} $A_{\theta }$ s'identifie isom\'{e}triquement 
\`{a} un sous espace de $A^{\theta }$.

\noindent\ \ \ \ \ D'apr\`{e}s \cite[theorem 4.2.2]{BL}, $A_{0}\cap A_{1}$
est toujours dense dans $A_{\theta },0<\theta <1.$ Si $A_{0}\cap A_{1}$ est
dense dans $A_{0}$ et $A_{1},$ on a ($A_{0}\cap A_{1})^{\ast }=A_{0}^{\ast
}+A_{1}^{\ast },$ $A_{0}^{\ast }\cap A_{1}^{\ast }=(A_{0}+A_{1})^{\ast }$ 
\cite[theorem 2.7.1]{BL}, on peut appliquer le th\'{e}or\`{e}me de reit\'{e}%
ration \cite[theorem 4.6.1]{BL} et ($A_{\theta })^{\ast }=(A_{0}^{\ast
},A_{1}^{\ast })^{\theta },$ $\theta \in \left] 0,1\right[ $ \cite[theorem
4.5.1]{BL}. On fait cette hypoth\`{e}se dans la suite.

Soient $\mu _{z}$ la mesure harmonique sur le bord de $S$ au point $z\in
S^{\circ }$ et $Q(z,.)$ le densit\'{e} de cette mesure par rapport \`{a} la
mesure de Lebesgue (sur les deux droites qui forment le bord de $S);$ notons 
$Q_{0}(z,\tau )=Q(z,i\tau )$ et $Q_{1}(z,\tau )=Q(z,1+i\tau ),$ $z\in
S^{\circ }$ et $\tau \in \mathbb{R}.$

\section{R\'{e}sultats}

Notons $B_{j}$ l'adh\'{e}rence de $A_{0}^{\ast }\cap A_{1}^{\ast }$ dans $%
A_{j}^{\ast },$ $j=0,1.$ Il est clair que $B_{0}\cap B_{1}=A_{0}^{\ast }\cap
A_{1}^{\ast },$ isom\'{e}triquement. D'apr\`{e}s \cite[Theorem 4.2.2 b)]{BL}
on a isom\'{e}triquement, pour $\theta \in \left] 0,1\right[ ,$

\begin{equation}
(B_{0},B_{1})_{\theta }=(A_{0}^{\ast },A_{1}^{\ast })_{\theta }\ .
\label{II}
\end{equation}

\noindent\ Comme $B_{0}\cap B_{1}$ est dense dans $B_{j},$ le dual de $%
B_{\theta }=(B_{0},B_{1})_{\theta }$ est $(B_{0}^{\ast },B_{1}^{\ast
})^{\theta }$ \cite[theorem 4.5.1]{BL} et, d'apr\`{e}s \cite[theorem 2.7.1]%
{BL},

\begin{equation*}
B_{0}^{\ast }+B_{1}^{\ast }=(B_{0}\cap B_{1})^{\ast }=(A_{0}^{\ast }\cap
A_{1}^{\ast })^{\ast }=(A_{0}+A_{1})^{\ast \ast }.
\end{equation*}

\noindent En particulier, $A_{0}+A_{1}$ s'identifie isom\'{e}triquement \`{a}
un sous espace ferm\'{e} de $B_{0}^{\ast }+B_{1}^{\ast }.$

\noindent Soit $i_{j}:B_{j}\rightarrow A_{j}^{\ast }$ l'application identit%
\'{e}; la restriction de son adjoint $i_{j}^{\ast }:A_{j}\rightarrow
B_{j}^{\ast },$ $j=0,1,$ est contractante.

\begin{lemma}
\label{R}Soit $R:$ $Q\mathcal{G}(A_{0},A_{1})\rightarrow Q\mathcal{G}%
(B_{0}^{\ast },B_{1}^{\ast }),$ l'application d\'{e}finie par $%
g(j+i.)\rightarrow $ $i_{j}^{\ast }(g(j+i.)),$ $j=0,1.$ $R$ \ est une
contraction et induit une contraction injective
\end{lemma}

\begin{equation*}
R^{\theta }:A^{\theta }\rightarrow (B_{0}^{\ast },B_{1}^{\ast })^{\theta
},\;\theta \in \left] 0,1\right[ .
\end{equation*}

\noindent \emph{D\'{e}monstration:} Il est clair que $R$ est une contraction
(non injective en g\'{e}n\'{e}ral). On identifie $A^{\theta }$ et $%
(B_{0}^{\ast },B_{1}^{\ast })^{\theta }$ \`{a} des quotients de $Q\mathcal{G}%
(A_{0},A_{1})$ et $Q\mathcal{G}(B_{0}^{\ast },B_{1}^{\ast })$
respectivement. Notant que $(R(g^{.}))^{\prime }(\theta )=R^{\theta
}(g^{\prime }(\theta )),$ $R$ induit une contraction $R^{\theta }$ sur ces
quotients. Notons que, pour $a\in A^{\theta },$ pour $b\in B_{0}\cap
B_{1}=A_{0}^{\ast }\cap A_{1}^{\ast }=(A_{0}+A_{1})^{\ast }$ (espace dense
dans $B_{\theta }),$

\begin{equation*}
\left\langle R^{\theta }(a),b\right\rangle =\left\langle a,b\right\rangle .
\end{equation*}%
Si $R^{\theta }(a)=0,$ alors $\left\langle a,b\right\rangle =0$ pour tout $%
b\in B_{0}\cap B_{1}=(A_{0}+A_{1})^{\ast },$ d'o\`{u} $a=0$ dans $%
A_{0}+A_{1},$ et dans $A^{\theta }$.$\blacksquare $

\begin{theorem}
\label{ci}Soient $(A_{0},A_{1})$ un couple d'interpolation complexe, $B_{j}$
l'adh\'{e}rence de $A_{0}^{\ast }\cap A_{1}^{\ast }$ dans $A_{j}^{\ast },$ $%
j=0,1,$ $\beta \in \left] 0,1\right[ ,$ $R^{\beta }$ d\'{e}finie comme
ci-dessus. Soit $Z^{\beta }$ l'adh\'{e}rence de $R^{\beta }(A^{\beta })$
dans $(B_{0}^{\ast },B_{1}^{\ast })^{\beta }.$ Supposons que $Z^{\beta }$
est un espace $W.C.G$ \cite[chap.VIII,p.251]{DU}$.$ Alors $A^{\theta
}=A_{\theta },$ pour tout $\theta \in \left] 0,1\right[ .$
\end{theorem}

\begin{remark}
\label{uu}Dans \cite{Da} on montre que si $A^{\beta }$ est un espace $W.C.G$
pour un $\beta \in \left] 0,1\right[ ,$ alors $A^{\theta }=A_{\theta },$
pour tout $\theta \in \left] 0,1\right[ .$ Remarquons que si $A^{\beta }$
est un espace $W.C.G,$ $Z^{\beta }$ est un espace $W.C.G.$
\end{remark}

\noindent La d\'{e}monstration du th\'{e}or\`{e}me \ref{ci}, n\'{e}cessite
les lemmes suivants.

\begin{lemma}
\label{yu}\textbf{\ }Pour tout $\theta \in \left] 0,1\right[ ,$ $R^{\theta }$
est une isom\'{e}trie: $A_{\theta }\rightarrow (B_{0}^{\ast },B_{1}^{\ast
})^{\theta }.$
\end{lemma}

\noindent \emph{D\'{e}monstration:} Comme $A_{\theta }$ s'identifie \`{a} un
sous espace de $A^{\theta }$ \cite{B}, $R^{\theta }$ est contractante: $%
A_{\theta }=(A_{0},A_{1})_{\theta }\rightarrow (B_{0}^{\ast },B_{1}^{\ast
})^{\theta }$ par le lemme \ref{R}.

\noindent Comme $A_{0}\cap A_{1}$ est dense dans $A_{\theta },$ il suffit de
montrer que $\left\Vert a\right\Vert _{A_{\theta }}\leq \left\Vert R^{\theta
}(a)\right\Vert _{(B_{0}^{\ast },B_{1}^{\ast })^{\theta }}$ lorsque $a\in
A_{0}\cap A_{1}.$

\noindent Soit $\varepsilon >0$; comme ($A_{\theta })^{\ast }=(A_{0}^{\ast
},A_{1}^{\ast })^{\theta },$ il existe $g\in \mathcal{G(}A_{0}^{\ast
},A_{1}^{\ast })$ tel que

\begin{equation}
\left\Vert a\right\Vert _{A_{\theta }}<\left\vert \left\langle a,g^{\prime
}(\theta )\right\rangle \right\vert +\varepsilon ,\text{ \qquad }\left\Vert
g^{.}\right\Vert _{Q\mathcal{G}(A_{0}^{\ast },A_{1}^{\ast })}\leq 1.
\label{io}
\end{equation}

\noindent Soient

\begin{equation*}
F_{n}(z)=-in\left[ g(z+i/n)-g(z)\right] ,\quad z\in S
\end{equation*}%
et $F_{n,\delta }(z)=e^{\delta z^{2}}F_{n}(z)$ $\ $pour $\delta >0.$ Comme $%
\left\vert F_{n}\right\vert $ est born\'{e}e sur le bord de $S,$ $\left\vert
F_{n,\delta }\right\vert $ tend vers $0$ \`{a} l'infini sur le bord, d'o\`{u}
$F_{n,\delta }\in \mathcal{F}(A_{0}^{\ast },A_{1}^{\ast }).$ Par d\'{e}%
finition 
\begin{equation*}
\left\Vert F_{n,\delta }(\theta )\right\Vert _{(A_{0}^{\ast },A_{1}^{\ast
})_{\theta }}\leq \left\Vert F_{n,\delta }\right\Vert _{\mathcal{F}%
(A_{0}^{\ast },A_{1}^{\ast })}\leq e^{\delta }\sup_{S}\left\vert
F_{n}(z)\right\vert \leq e^{\delta }\left\Vert g^{.}\right\Vert _{\mathcal{QG%
}(A_{0}^{\ast },A_{1}^{\ast })}\leq e^{\delta }.
\end{equation*}

\noindent D'o\`{u}, pour tout $n,$ par (\ref{II}),

\begin{equation*}
\left\Vert F_{n}(\theta )\right\Vert _{(B_{0},B_{1})_{\theta }}=\left\Vert
e^{-\delta \theta ^{2}}F_{n,\delta }(\theta )\right\Vert _{(A_{0}^{\ast
},A_{1}^{\ast })_{\theta }}=\lim_{\delta \rightarrow 0}\left\Vert e^{-\delta
\theta ^{2}}F_{n,\delta }(\theta )\right\Vert _{(A_{0}^{\ast },A_{1}^{\ast
})_{\theta }}\leq 1.
\end{equation*}

\noindent Comme $g$ est holomorphe \`{a} valeurs dans $A_{0}^{\ast
}+A_{1}^{\ast }=(A_{0}\cap A_{1})^{\ast }$, $\left\langle a,F_{n}(\theta
)\right\rangle \underset{n\rightarrow \infty }{\rightarrow }\left\langle
a,g^{\prime }(\theta )\right\rangle .$ Il existe $n_{0}$ assez grand tel
que, d'apr\`{e}s (\ref{io}),

\begin{equation*}
-2\varepsilon +\left\Vert a\right\Vert _{A_{\theta }}<\left\vert
\left\langle a,F_{n_{0}}(\theta )\right\rangle \right\vert \leq \left\Vert
R^{\theta }(a)\right\Vert _{(B_{0}^{\ast },B_{1}^{\ast })^{\theta
}}\left\Vert F_{n_{0}}(\theta )\right\Vert _{(B_{0},B_{1})_{\theta }}\leq
\left\Vert R^{\theta }(a)\right\Vert _{(B_{0}^{\ast },B_{1}^{\ast })^{\theta
}},
\end{equation*}

\noindent d$^{\prime }$o\`{u} l'in\'{e}galit\'{e} cherch\'{e}e lorsque $%
\varepsilon \rightarrow 0.\blacksquare $

\begin{lemma}
\label{bn} Soient $g\in \mathcal{G(}\overline{A})$, $\theta \in \left] 0,1%
\right[ .$ L'application$:\tau \rightarrow R^{\theta }(g^{\prime }(\theta
+i\tau ))$ est born\'{e}e: $\mathbb{R}\rightarrow (B_{0}^{\ast },B_{1}^{\ast
})^{\theta }.$ Pour tout $c\in (B_{0}^{\ast },B_{1}^{\ast })^{\theta },$
l'application$:\tau \rightarrow \left\Vert c+R^{\theta }(g^{\prime }(\theta
+i\tau ))\right\Vert _{(B_{0}^{\ast },B_{1}^{\ast })^{\theta }}$ est s.c.i
sur $\ \mathbb{R}.$
\end{lemma}

\noindent \emph{D\'{e}monstration:} Par d\'{e}finition de $A^{\theta }$, $%
g^{\prime }(\theta )\in A^{\theta }$ ; par le lemme \ref{R}

\begin{equation*}
\left\Vert R^{\theta }(g^{\prime }(\theta ))\right\Vert _{(B_{0}^{\ast
},B_{1}^{\ast })^{\theta }}\leq \left\Vert g^{\prime }(\theta )\right\Vert
_{A^{\theta }}\leq \left\Vert g^{.}\right\Vert _{Q\mathcal{G}(\overline{A})}.
\end{equation*}

\noindent La fonction $g_{i\tau }$ d\'{e}finie par $g_{i\tau }(z)=g(z+i\tau
),$ $z\in S,$ $\tau \in \mathbb{R},$ v\'{e}rifie $\left\Vert g_{i\tau
}^{.}\right\Vert _{Q\mathcal{G}(\overline{A})}=\left\Vert g^{.}\right\Vert
_{Q\mathcal{G}(\overline{A})},$ donc $\left\Vert R^{\theta }(g_{i\tau
}^{\prime }(\theta ))\right\Vert _{(B_{0}^{\ast },B_{1}^{\ast })^{\theta
}}\leq \left\Vert g^{.}\right\Vert _{Q\mathcal{G}(\overline{A})}.$

\noindent D'apr\`{e}s (\ref{II}), et comme $B_{0}\cap B_{1}=A_{0}^{\ast
}\cap A_{1}^{\ast }$ est dense dans $B_{\theta },$ on a

\begin{eqnarray*}
&&\left\Vert c+R^{\theta }(g^{\prime }(\theta +i\tau ))\right\Vert
_{(B_{0}^{\ast },B_{1}^{\ast })^{\theta }} \\
&=&sup\left\{ \left\vert \left\langle b,c+R^{\theta }(g^{\prime }(\theta
+i\tau ))\right\rangle \right\vert ;\text{ }\left\Vert b\right\Vert
_{(B_{0},B_{1})_{\theta }}\leq 1\right\} \\
&=&sup\left\{ \left\vert \left\langle a^{\ast },c+g^{\prime }(\theta +i\tau
)\right\rangle \right\vert ;\text{ }a^{\ast }\in A_{0}^{\ast }\cap
A_{1}^{\ast },\quad \left\Vert a^{\ast }\right\Vert _{(A_{0}^{\ast
},A_{1}^{\ast })_{\theta }}\leq 1\right\} .
\end{eqnarray*}%
Comme $g$ est holomorphe \`{a} valeurs dans $A_{0}+A_{1},$ pour tout $%
a^{\ast }\in A_{0}^{\ast }\cap A_{1}^{\ast }=(A_{0}+A_{1})^{\ast },$ les
applications $\tau $ $\rightarrow \left\vert \left\langle a^{\ast
},a+g^{\prime }(\theta +i\tau )\right\rangle \right\vert $ sont continues
sur $\mathbb{R}$. Leur supremum est donc s.c.i.$.\blacksquare $

Par un argument analogue \`{a} celui de \cite{Da}, on montre

\begin{lemma}
\label{k}Soient $\overline{C}=(C_{0},C_{1})$ un couple d'interpolation, $%
\beta \in \left] 0,1\right[ $, $Z^{\beta }$ un sous-espace $W.C.G$ $,$ $g\in 
\mathcal{G(}\overline{C})$. Alors l'application$:\tau \in \mathbb{R}%
\rightarrow R^{\beta }(g^{\prime }((\beta +i\tau ))$ est p.s \'{e}gale \`{a}
une fonction fortement mesurable: $\mathbb{R}\rightarrow Z^{\beta }\subset
C_{\beta }^{\ast \ast }.$
\end{lemma}

\begin{lemma}
\label{mesurable}\textbf{\ }Soient $g\in \mathcal{G(}\overline{A})$, $\phi
_{\theta }(t)=g^{\prime }(\theta +it),t\in \mathbb{R}$.

i) Si $\phi _{\theta }$ est \`{a} valeurs dans un sous espace ferm\'{e} s%
\'{e}parable $Z$ de $A_{\theta },$ elle est fortement mesurable: $\mathbb{R}%
\rightarrow A_{\theta }.$

ii) Si $\phi _{\theta }$ est \`{a} valeurs dans un sous espace ferm\'{e} s%
\'{e}parable $Z$ de $A^{\theta },$ elle est fortement mesurable: $\mathbb{R}%
\rightarrow $ $A^{\theta }.$
\end{lemma}

\noindent Dans la suite on utilise seulement i), dans la preuve du lemme \ref%
{lem} d). On donne deux preuves de i) (noter que ii) implique i)).

\noindent \emph{Preuve:} i) D'apr\`{e}s le lemme \ref{yu}, $Z$ est un sous
espace ferm\'{e} de $(B_{0}^{\ast },B_{1}^{\ast })^{\theta }$ et, d'apr\`{e}%
s le lemme \ref{bn}, l'application $t\rightarrow \left\Vert \phi _{\theta
}(t)-c\right\Vert _{Z}$ est s.c.i. pour tout $c\in Z.$ L'image r\'{e}%
ciproque par $\phi _{\theta }$ de toute boule ouverte de $Z$ est donc un bor%
\'{e}lien. Comme $Z$ est s\'{e}parable, tout ouvert de $Z$ est r\'{e}union d%
\'{e}nombrable de boules, donc $\phi _{\theta }$ est bien mesurable \`{a}
valeurs dans $Z$.$\blacksquare $

ii) Soient $J$ l'injection canonique: $Z\rightarrow A_{0}+A_{1},$ et $Y$ \
l'adh\'{e}rence de $J(Z)$ dans $A_{0}+A_{1}.$ Comme $Z$ et $Y$ sont des
espaces polonais, comme $J$ est continue, $J^{-1}$ est bor\'{e}lienne: $%
J(Z)\rightarrow Z,$ voir par exemple \cite{A}. Comme $J\circ \phi _{\theta }:%
\mathbb{R}\rightarrow A_{0}+A_{1}$ est continue et \`{a} valeurs dans $J(Z),$
comme $\phi _{\theta }=J^{-1}\circ (J\circ \phi _{\theta }),$ alors $\phi
_{\theta }$ est bor\'{e}lienne: $\mathbb{R}\rightarrow Z.\blacksquare $

\begin{lemma}
\label{lem} Soient $g\in \mathcal{G(}\overline{A})$, $\beta \in \left] 0,1%
\right[ $ et $\phi _{\beta }=g^{\prime }(\beta +i.).$

a) On suppose que $R^{\beta }\circ \phi _{\beta }$ est p.s. \'{e}gale \`{a}
une fonction fortement mesurable: $\mathbb{R}\rightarrow (B_{0}^{\ast
},B_{1}^{\ast })^{\beta }.$ Alors $\phi _{\beta }$ est p.s. \`{a} valeurs
dans $A_{\beta }.$

On suppose d\'{e}sormais que $\phi _{\beta }$ est p.s. \'{e}gale \`{a} une
fonction fortement mesurable: $\mathbb{R}\rightarrow A_{\beta }.$ Alors

b) pour $\theta \neq \beta $ $\ g^{\prime }(\theta )\in A_{\theta }.$

c) pour tout $\theta \neq \beta ,$ $\phi _{\theta }$ est \`{a} valeurs dans
un sous espace s\'{e}parable de $A_{\theta }.$

d) $g^{\prime }(\beta )\in A_{\beta }.$
\end{lemma}

\noindent On a not\'{e} $R^{\beta }\circ \phi _{\beta }$ la fonction: $%
t\rightarrow R^{\beta }(g_{it}^{\prime }(\beta )).$

\noindent \emph{Preuve:} a) \emph{\'{e}tape 1: }Posons

\begin{equation*}
g_{1}=g-g(0)-\alpha _{0}
\end{equation*}%
o\`{u} $g(1)-g(0)=\alpha _{0}+\alpha _{1}$ ($\alpha _{j}\in A_{j}$, $j=0,1),$
avec

\begin{equation*}
\left\Vert g(1)-g(0)\right\Vert _{A_{0}+A_{1}}=\left\Vert \alpha
_{0}\right\Vert _{A_{0}}+\left\Vert \alpha _{1}\right\Vert _{A_{1}}.
\end{equation*}

Comme $g$ est holomorphe \`{a} l'int\'{e}rieur de $S,$ pour tous $t\in 
\mathbb{R},$ $h>0,$ $\theta \in ]0,1[,$ on a, dans $A_{0}+A_{1},$

\begin{equation}
g(\theta +i(t+h))-g(\theta +it)=\int_{t}^{t+h}g^{\prime }(\theta +i\tau
)d\tau  \label{4 bis}
\end{equation}

\noindent

\noindent D'apr\`{e}s l'in\'{e}galit\'{e} des accroissements finis et (\ref%
{2})

\begin{equation*}
\left\Vert g(1)-g(0)\right\Vert _{A_{0}+A_{1}}\leq \left\Vert g^{\cdot
}\right\Vert _{Q\mathcal{G}(\overline{A})}.
\end{equation*}%
Alors $g_{1}:S\rightarrow $ $A_{0}+A_{1}$ est continue sur $S$ et holomorphe 
\`{a} l'int\'{e}rieur de $S.$ Comme $g\in \mathcal{G(}\overline{A}),$ pour
tout $\tau \in \mathbb{R}$ et $j\in \{0,1\},$ on a

\begin{equation*}
\left\Vert g_{1}(j+i\tau )\right\Vert _{A_{j}}\leq \left\Vert g(j+i\tau
)-g(j)\right\Vert _{A_{j}}+\left\Vert \alpha _{j}\right\Vert _{A_{j}}\leq
(1+\left\vert \tau \right\vert )\left\Vert g^{\cdot }\right\Vert _{Q\mathcal{%
G}(\overline{A})}.
\end{equation*}

\noindent L'application $z\rightarrow G_{\varepsilon }(z)=e^{\varepsilon
z^{2}}g_{1}(z)$ est donc dans $\mathcal{F}\overline{(A})$ pour tout $%
\varepsilon >0.$ En particulier, pour tout $t\in \mathbb{R},$ $%
G_{\varepsilon }(\theta +it)\in A_{\theta },$ donc $g_{1}(\theta +it)\in
A_{\theta }$. D'o\`{u}

\begin{equation*}
g_{1}(\theta +i(t+h))-g_{1}(\theta +it)=g(\theta +i(t+h))-g(\theta +it)\in
A_{\theta }.
\end{equation*}

\noindent Alors, d'apr\`{e}s (\ref{4 bis}), $\int_{t}^{t+h}g^{\prime
}(\theta +i\tau )d\tau $ est dans $A_{\theta },$ pour $t$ et $h$ r\'{e}els.

\emph{\'{e}tape 2: } Par hypoth\`{e}se $R^{\beta }\circ \phi _{\beta }$ est
p.s. \'{e}gale \`{a} une fonction fortement mesurable: $\mathbb{R}%
\rightarrow Z^{\beta },$ o\`{u} $Z^{\beta }$ est l'adh\'{e}rence de $%
A^{\beta }$ dans $(B_{0}^{\ast },B_{1}^{\ast })^{\beta }.$ Le th\'{e}or\`{e}%
me de differentiabilit\'{e} de Lebesgue \cite[chap.II th.9 p. 48]{DU} entra%
\^{\i}ne que, p.s., on a dans $Z^{\beta }$ \ l'\'{e}galit\'{e}

\begin{equation}
iR^{\beta }\circ \phi _{\beta }(it)=\lim_{h\rightarrow 0}\frac{1}{h}%
\int_{t}^{t+h}R^{\beta }\circ \phi _{\beta }(i\tau )d\tau
=\lim_{h\rightarrow 0}R^{\beta }(\frac{1}{h}\int_{t}^{t+h}g^{\prime }(\beta
+i\tau )d\tau ),  \label{Leb}
\end{equation}

\noindent o\`{u} $h$ est r\'{e}el. D'apr\`{e}s la fin de l'\'{e}tape 1
appliqu\'{e}e en $\beta $ et le lemme \ref{yu}, cette limite dans $Z^{\beta
} $ est en fait une limite dans $A_{\beta },$ c\`{a}d p.s. $g^{\prime
}(\beta +i.)\in $ $A_{\beta }.$

b) On suppose d'abord $\theta >\beta .$

\noindent \emph{\'{e}tape 1: } Soit

\begin{equation*}
V(z)=g_{1}(\beta +(1-\beta )z),\;z\in S.
\end{equation*}

\noindent Cette fonction est \`{a} valeurs dans $A_{0}+A_{1}$ est holomorphe 
\`{a} l'int\'{e}rieur de $S$ et continue sur $S,$ donc s'exprime \`{a}
l'aide de la mesure harmonique sur le bord de $S.$\textbf{\ }Pour\textbf{\ }v%
\'{e}rifier\textbf{\ }que $V,$ vue comme fonction \`{a} valeurs dans $%
A_{\beta }+A_{1},$ est holomorphe \`{a} l'int\'{e}rieur de $S$ et continue
sur $S,$ il suffira donc de voir que $V$ est continue sur l'axe imaginaire, 
\`{a} valeurs dans $A_{\beta }.$

On va montrer que $V\in \mathcal{G}(A_{\beta },A_{1})$ avec une norme $\leq
(1-\beta )\left\Vert g^{\cdot }\right\Vert _{Q\mathcal{G}(\overline{A})}.$
L'in\'{e}galit\'{e} correspondante sur la droite Re $z=1$ est \'{e}vidente.
Pour la v\'{e}rifier sur l'axe imaginaire, posons, pour $\tau ,\tau ^{\prime
}$ r\'{e}els fix\'{e}s,

\begin{equation*}
F_{\tau ,\tau ^{\prime }}(\xi )=\frac{g(\xi +i(1-\beta )\tau )-g(\xi
+i(1-\beta )\tau ^{\prime })}{\tau -\tau ^{\prime }},\;\xi \in S,
\end{equation*}

\noindent d'o\`{u} $F_{\tau ,\tau ^{\prime }}(\beta )=\frac{V(i\tau
)-V(i\tau ^{\prime })}{\tau -\tau ^{\prime }}$, et $F_{\tau ,\tau ^{\prime
}}(1)=\frac{V(1+i\tau )-V(1+i\tau ^{\prime })}{\tau -\tau ^{\prime }}$. Pour
tout $t\in \mathbb{R},$ on a

\begin{equation*}
\left\Vert F_{\tau ,\tau ^{\prime }}(j+it)\right\Vert _{A_{j}}\leq (1-\beta
)\left\Vert g^{\cdot }\right\Vert _{Q\mathcal{G}(\overline{A})},j\in \{0,1\}.
\end{equation*}

\noindent Comme dans l'\'{e}tape 1 de a), pour tout $\varepsilon >0,$
l'application $\xi \rightarrow H_{\varepsilon ,\tau ,\tau ^{\prime }}(\xi
)=e^{\varepsilon \xi ^{2}}F_{\tau ,\tau ^{\prime }}(\xi )$ v\'{e}rifie

\begin{equation*}
\left\Vert H_{\varepsilon ,\tau ,\tau ^{\prime }}\right\Vert _{\mathcal{F}%
\overline{(A})}\leq e^{\varepsilon }(1-\beta )\left\Vert g^{\cdot
}\right\Vert _{Q\mathcal{G}(\overline{A})},
\end{equation*}

\noindent d'o\`{u}

\begin{equation*}
\left\Vert F_{\tau ,\tau ^{\prime }}(\beta )\right\Vert _{A_{\beta }}\leq
(1-\beta )\left\Vert g^{\cdot }\right\Vert _{Q\mathcal{G}(\overline{A})}.
\end{equation*}

\noindent On a donc, pour tous $\tau ,\tau ^{\prime }$ r\'{e}els,

\begin{equation*}
\left\Vert V(i\tau )-V(i\tau ^{\prime })\right\Vert _{A_{\beta }}\leq
\left\vert \tau -\tau ^{\prime }\right\vert (1-\beta )\left\Vert g^{\cdot
}\right\Vert _{Q\mathcal{G}(\overline{A})},
\end{equation*}

\noindent ce qui prouve la continuite de $V$ sur l'axe imaginaire, \`{a}
valeurs dans $A_{\beta },$ et l'assertion annonc\'{e}e.

\emph{\'{e}tape 2: } d'apr\`{e}s la preuve de a), pour $h$ r\'{e}el, p.s.

\begin{equation*}
\lim_{h\rightarrow 0}(V(i(\tau +h))-V(i\tau ))/h=(1-\beta )g^{\prime }(\beta
+(1-\beta )i\tau )\quad dans\quad A_{\beta }.
\end{equation*}

\noindent D'apr\`{e}s \cite[lemma 4.3.3]{BL}, on a alors

\begin{equation*}
V^{\prime }(\eta )\in (A_{\beta },A_{1})_{\eta },\eta \in \left] 0,1\right[ .
\end{equation*}

\emph{\'{e}tape 3:} Choisissons $\eta $ tel que $\theta =(1-\eta )\beta
+\eta $. D'apr\`{e}s le th\'{e}or\`{e}me de r\'{e}it\'{e}ration \cite[%
theorem 4.6.1]{BL}, $(A_{\beta },A_{1})_{\eta }=A_{\theta }$, donc

\begin{equation*}
V^{\prime }(\eta )=(1-\beta )g^{\prime }(\theta )\in A_{\theta },
\end{equation*}

\noindent ce qui ach\`{e}ve la preuve lorsque $\beta <\theta .$

Si $0<\theta <\beta $ le raisonnement est analogue, en rempla\c{c}ant $V$
par $W(z)=g_{1}(\beta z)\in \mathcal{G}(A_{0},A_{\beta })$, telle que $%
\lim_{h\rightarrow 0}(W(1+i(\tau +h))-W(1+i\tau ))/h$ existe dans $A_{\beta
},$ pour presque tout $\tau ,$ avec $h$ r\'{e}el.

c) Soit $A_{0}^{\prime }\subset A_{0}$ \ le sous espace ferm\'{e} s\'{e}%
parable engendr\'{e} par $\{g_{1}(it)$, $t\in \mathbb{R\}}$. Comme $g_{1}$
est continue sur $S,$ $A_{0}^{\prime }$ est s\'{e}parable, ainsi que $%
(A_{0}^{\prime },A_{1})_{\beta }$ et son adh\'{e}rence $Y$ dans $A_{\beta }.$
Par l'\'{e}tape 2 de b) appliqu\'{e}e au couple $(A_{0}^{\prime },A_{1})$, $%
g^{\prime }(\beta +it)$ est p.s. dans $(A_{0}^{\prime },A_{1})_{\beta },$
donc p.s. dans $Y,$ ce qui r\`{e}gle le cas $\theta =\beta .$

Pour le cas $\beta <\theta ,$ rempla\c{c}ons la fonction $V$ de l'\'{e}tape
1 de b) par $V_{t}(z)=V(z+it)$, avec $t$ fix\'{e} r\'{e}el. Comme en b), $%
V_{t}\in \mathcal{G}(Y,A_{1}),$ $V_{t}^{\prime }(\eta )\in (Y,A_{1})_{\eta
},\eta \in \left] 0,1\right[ $ et $(Y,A_{1})_{\eta }$ est s\'{e}parable.
Soit $\eta $ d\'{e}fini comme dans l'\'{e}tape 3 de b)$.$ Comme ci-dessus, $%
V_{t}^{\prime }(\eta )=(1-\beta )g^{\prime }(\theta +i(1-\beta )t).$ Soit $%
Z_{\theta }$ l'adh\'{e}rence de $(Y,A_{1})_{\eta }$ dans $(A_{\beta
},A_{1})_{\eta }=A_{\theta }$; $Z_{\theta }$ est donc s\'{e}parable et $\phi
_{\theta }=g^{\prime }(\theta +i.)$ est \`{a} valeurs dans $Z_{\theta }.$

On raisonne de fa\c{c}on analogue si $0<\theta <\beta $ en consid\'{e}rant $%
W_{t}(z)=W(z+it):$ $W_{t}$ est dans $\mathcal{G}(A_{0}^{\prime },Y).$

d) Soit $\theta >\beta .$ Par c) et le lemme \ref{mesurable} i), $\phi
_{\theta }$ est fortement mesurable \`{a} valeurs dans $A_{\theta }.$ Alors
b) appliqu\'{e} en \'{e}changeant les r\^{o}les de $\beta $ et $\theta $
donne $g^{\prime }(\beta )\in A_{\beta }.\blacksquare $

\noindent \emph{D\'{e}monstration du th\'{e}or\`{e}me 2:} Soient $a\in
A^{\beta }$ et $g\in \mathcal{G}(\overline{A})$ tels que $a=g^{\prime
}(\beta ).$ D'apr\`{e}s le lemme \ref{R}, l'application $R^{\beta }\circ
\phi _{\beta }:\tau \in \mathbb{R}\rightarrow R^{\beta }(g^{\prime }(\beta
+i\tau ))$ est \`{a} valeurs dans $Z^{\beta }\subset (B_{0}^{\ast
},B_{1}^{\ast })^{\beta }.$ Gr\^{a}ce \`{a} l'hypoth\`{e}se sur $Z^{\beta },$
on peut appliquer le lemme \ref{k}, donc $R^{\beta }\circ \phi _{\beta }$
est p.s. \'{e}gale \`{a} une fonction fortement mesurable \`{a} valeurs dans 
$(B_{0}^{\ast },B_{1}^{\ast })^{\beta },$ et $\phi _{\beta }$ est p.s. \'{e}%
gale \`{a} une fonction fortement mesurable \`{a} valeurs dans $A_{\beta }$
par le lemme \ref{lem} a). D'apr\`{e}s le lemme \ref{lem} b) $g^{\prime
}(\theta )\in A_{\theta }$ pour tout $\theta \neq \beta .$ Il en r\'{e}sulte
que $A^{\theta }=A_{\theta },$ pour tout $\theta \neq \beta .$ Enfin par le
lemme \ref{lem} d) $g^{\prime }(\beta )=a\in A_{\beta },$ d'o\`{u} $A^{\beta
}=A_{\beta }.\blacksquare $

\begin{proposition}
\label{jj}Soit $(A_{0},A_{1})$ un couple d'interpolation; supposons qu'il
existe un $\beta $ tel que
\end{proposition}

1)$(A_{0}^{\ast },A_{1}^{\ast })_{\beta }$ est s\'{e}parable.

2)$(A_{0}^{\ast },A_{1}^{\ast })_{\beta }$ ne contient pas $\ell ^{1}$
isomorphiquement.

3)$A_{\beta }$ est un espace compl\'{e}ment\'{e} de $\left[ (A_{0}^{\ast
},A_{1}^{\ast })_{\beta }\right] ^{\ast }$ (d'apr\`{e}s lemme \ref{yu}, $%
A_{\beta }$ est un sous-espace isom\'{e}trique de $\left[ (A_{0}^{\ast
},A_{1}^{\ast })_{\beta }\right] ^{\ast }).$

Alors $A^{\theta }=A_{\theta },$ pour tout $\theta \in \left] 0,1\right[ .$

\begin{lemma}
\label{nn}Supposons qu'il existe un $\beta \in \left] 0,1\right[ $ tel que $%
A_{\beta }$ a la propri\'{e}t\'{e} de Radon-Nikodym \cite{E}. Alors $%
A^{\theta }=A_{\theta },$ pour tout $\theta \in \left] 0,1\right[ .$
\end{lemma}

\bigskip \emph{Preuve:}Soit $\theta ,\gamma ,\eta \in 0,1$ $(\gamma >\beta )$
tel que $\theta =(\eta -1)\beta +\eta \gamma ;$ montrons d'abord que $%
A^{\theta }\subset (A_{\beta },A_{\gamma })^{\eta }.$

Pour cela, soient $a\in $ $(A_{0},A_{1})^{\theta }$ et $g\in \mathcal{G}%
(A_{0},A)$ tels que $g^{\prime }(\theta )=a.$

Posons

\begin{equation*}
g_{1}=g-g(0)-\alpha _{0}
\end{equation*}%
o\`{u} $g(1)-g(0)=\alpha _{0}+\alpha _{1}$ ($\alpha _{j}\in A_{j}$, $j=0,1),$
avec

\begin{equation*}
\left\Vert g(1)-g(0)\right\Vert _{A_{0}+A_{1}}=\left\Vert \alpha
_{0}\right\Vert _{A_{0}}+\left\Vert \alpha _{1}\right\Vert _{A_{1}}.
\end{equation*}

Soit $\alpha \in \left] 0,1\right[ $ tel que $\alpha +\beta =\gamma $

\begin{equation*}
V(z)=g_{1}(\alpha z+\beta ),\;z\in S.
\end{equation*}

D'apr\`{e}s (a) \'{e}tape 1 du lemme \ref{lem}, $V\in \mathcal{G}(A_{\beta
},A_{\gamma }),$ donc $V^{\prime }(\eta )=\alpha g^{\prime }(\theta )\in
(A_{\beta },A_{g})^{\eta },$ d'o\`{u} $A^{\theta }\subset (A_{\beta
},A_{\gamma })^{\eta }.$

Supposons maintenant que $A_{\beta }$ a la propri\'{e}t\'{e} de Radon-Nikodym%
$;$ d'apr\`{e}s le r\'{e}sultat de \cite[Corol.4.5.2]{BL} on a $(A_{\beta
},A_{\gamma })^{\eta }=(A_{\beta },A_{\gamma })_{\eta }.$ D'autre part par
le th\'{e}or\`{e}me de reit\'{e}ration \cite[th.4.6.1]{BL} on a $(A_{\beta
},A_{\gamma })_{\eta }=A_{\theta },$ ce qui implique $A^{\theta }=A_{\theta
}.\blacksquare $

D\'{e}monstration de la proposition \ref{jj}: D'apr\`{e}s le r\'{e}sultat 
\cite{GH-GOD-MAU-SCH} et les conditions (1), (2) et (3), $A_{\beta }$ a la
propri\'{e}t\'{e} de Radon-Nikodym. Pour conclure le th\'{e}or\`{e}me, il
suffit d'appliquer le lemme \ref{nn}.$\blacksquare $

Par un argument analogue \`{a} celui de la proposition \ref{jj} on montre,

\begin{proposition}
\label{tt}Soit $(A_{0},A_{1})$ un couple d'interpolation; supposons qu'il
existe un $\beta \in \left] 0,1\right[ $ tel que
\end{proposition}

1)$(A_{0},A_{1})_{\beta }$ est s\'{e}parable.

2)$(A_{0},A_{1})_{\beta }$ ne contient pas $\ell ^{1}$ isomorphiquement.

3)$(A_{0}^{\ast },A_{1}^{\ast })_{\beta }$ est un espace compl\'{e}ment\'{e}
de $(A_{0}^{\ast },A_{1}^{\ast })^{\beta }.$

Alors $(A_{0}^{\ast },A_{1}^{\ast })^{\theta }=$ $(A_{0}^{\ast },A_{1}^{\ast
})_{\theta },$ pour tout $\theta \in \left] 0,1\right[ .$

$B_{1}(\mathbb{R})$ sera not\'{e} l'espace des fonctions:$\mathbb{%
R\rightarrow R}$ de premi\`{e}re classe de Baire.

Pour tout espace de Banach $X$, notons $B_{X}$ la boule unit\'{e} ferm\'{e}e
de $X.$

\begin{proposition}
\bigskip \label{gf}Soient $(A_{0},A_{1})$ un couple d'interpolation et $%
\beta \in \left] 0,1\right[ ;$ suppsons que \ \ \ \ \ \ \ \ \ \ \ \ \ \ \ \
\ \ \ \ \ \ \ \ \ \ \ \ \ \ \ \ \ \ \ \ \ \ \ \ \ \ \ \ \ \ \ \ \ \ \ \ \ \
\ \ \ \ \ \ \ \ \ \ \ \ \ \ \ \ \ \ \ \ \ \ \ \ \ \ \ \ \ \ \ \ \ \ \ \ \ \
\ \ \ \ \ \ \ \ \ \ \ \ \ \ \ \ \ \ \ \ \ \ \ \ \ \ \ \ \ \ \ \ \ \ \ \ \ \
\ \ \ \ \ \ \ \ \ \ \ \ \ \ \ \ \ \ \ \ \ \ \ \ \ \ \ \ \ \ \ \ \ \ \ \ \ \
\ \ \ \ \ \ \ \ \ \ \ \ \ \ \ \ \ \ \ \ \ \ \ \ \ \ \ \ \ \ \ \ \ \ \ \ \ \
\ \ \ \ \ \ \ \ \ \ \ \ \ \ \ \ \ \ \ \ \ \ \ \ \ \ \ \ \ \ \ \ \ \ \ \ \ \
\ \ \ \ \ \ \ 
\end{proposition}

1)$A_{\beta }$ est un espace faiblement de Lindel\"{o}f.

2)Il existe une projection continue $P:Z^{\beta }\rightarrow A_{\beta }.$

3)Pour tout $a^{\ast }\in (Z^{\beta })^{\ast },$ l'application$:\tau \in 
\mathbb{R}\rightarrow (R^{\beta }\circ \phi _{\beta }(\tau ),a^{\ast })\in 
\mathbb{C}$ est de premi\`{e}re classe de Baire.

Alors $A_{\theta }=A^{\theta },$ pour tout $\theta \in \left] 0,1\right[ .$

\emph{Preuve:}D'apr\`{e}s le lemme \ref{lem}, il suffit de montrer que
l'application $\tau \rightarrow R^{\beta }\circ g^{\prime }(\beta +i\tau )$
est \'{e}gale presque-partout \`{a} une fonction mesurable dans $Z^{\beta }.$

Soit $g\in \mathcal{G}(A_{0},A_{1});$ notons $F$ le sous-espace ferm\'{e}
engendr\'{e} par $\left\{ R^{\beta }\circ g^{\prime }(\beta +i\tau );\text{ }%
\tau \in \mathbb{R}\right\} $ dans $Z^{\beta }.$

\emph{\'{e}tape 1:}La boule unit\'{e} de $F^{\ast }$ est un compact de
Rosenthal \cite{GOD} pour la topologie pr\'{e}faible.

Soit $x^{\ast }\in B_{F^{\ast }};$ on d\'{e}finit $\sigma _{x^{\ast }}:%
\mathbb{R}\rightarrow \mathbb{C},$ par $\sigma _{x^{\ast }}(\tau )=(R^{\beta
}\circ g^{\prime }(\beta +i\tau ),x^{\ast }),$ $\tau \in \mathbb{R}.$ Consid%
\'{e}rons $U:B_{F^{\ast }}\rightarrow B_{1}(\mathbb{R})$, l'application d%
\'{e}finie par $U(x^{\ast })=\sigma _{x^{\ast }}.$ $U$ est continue
injective, donc $B_{F^{\ast }}$ est un compact de Rosenthal \cite{GOD} (on
peut supposer que $\sigma _{x^{\ast }}$ est \`{a} valeurs dans $\mathbb{R)}.$

\emph{\'{e}tape 2:}La boule unit\'{e} de $(A_{0}^{\ast },A_{1}^{\ast
})_{\beta }$ est pr\'{e}faiblement dense dans la boule unit\'{e} de $%
(Z^{\beta })^{\ast }.$

Soit $a^{\ast }\in B_{(Z^{\beta })^{\ast }};$ d'apr\`{e}s le th\'{e}or\`{e}%
me de Hahn-Banach, il existe $a^{\ast \ast }\in B_{(A_{0}^{\ast
},A_{1}^{\ast })_{\beta }^{\ast \ast }}$ qui prolonge $a^{\ast },$ pour
conclure, il suffit de remarquer que la boule unit\'{e} de $(A_{0}^{\ast
},A_{1}^{\ast })_{\beta }$ est pr\'{e}faiblement dense dans $B_{(A_{0}^{\ast
},A_{1}^{\ast })_{\beta }^{\ast \ast }}.$

\emph{\'{e}tape 3:} Il existe une fonction $\phi _{1}$ fortement mesurable 
\`{a} valeurs dans $A_{\beta }$ telle que $(P(R^{\beta }\circ \phi (\tau
),x^{\ast })=(\phi _{1}(\tau ),x^{\ast })$ pour presque-tout $\tau \in 
\mathbb{R}$ et tout \ $x^{\ast }\in A_{\beta }^{\ast }$

D'apr\`{e}s l'hypoth\`{e}se (3) $P(R^{\beta }\circ \phi _{\beta }(.))$ est
faiblement mesurable, l'espace $A_{\beta }$ est faiblement de Lindel\"{o}f,
donc il est mesure compact \cite{MOR1} par cons\'{e}quent il existe une
fonction $\phi _{1}$ fortement mesurable \`{a} valeurs dans $A_{\beta }$
telle que pour presque-tout $\tau \in \mathbb{R}$ et tout\ $x^{\ast }\in
A_{\beta }^{\ast }$%
\begin{equation}
(P(R^{\beta }\circ \phi _{\beta }(\tau ),x^{\ast })=(\phi _{1}(\tau
),x^{\ast })\text{ \cite{MOR1}-\cite{MOR2}.}  \label{p}
\end{equation}

\emph{\'{e}tape 4:}Pour tout $x^{\ast }\in (Z_{\beta })^{\ast }$ et tout $%
n\in \mathbb{N}^{\ast }$ $-i(R^{\beta }\circ g(\beta +i\tau +i/n)-R^{\beta
}\circ g(\beta +i\tau ),x^{\ast })=\left[ \dint\limits_{\tau }^{\tau +1/n}(P%
\left[ R^{\beta }\circ \phi _{\beta }(t)\right] ,x^{\ast })dt\right] .$

(remarquons d'apr\`{e}s le lemme \ref{yu} que $R^{\beta }\circ g(\beta
+i\tau +i/n),R^{\beta }\circ g(\beta +i\tau )\in A_{\beta },$ pour tout $%
\tau \in \mathbb{R}$ et tout $n\in \mathbb{N}^{\ast }).$

\bigskip Soit $x^{\ast }\in (Z_{\beta })^{\ast }$ tel que $P^{\ast }x^{\ast
}\in B_{(Z^{\beta })^{\ast }};$ d'apr\`{e}s l'\'{e}tape $2$, $P^{\ast
}x^{\ast }$ est $\sigma ((Z^{\beta })^{\ast },Z^{\beta })$ adh\'{e}rent \`{a}
la boule unit\'{e} de $(A_{0}^{\ast },A_{1}^{\ast })_{\beta },$ par cons\'{e}%
quent il existe une suite g\'{e}n\'{e}ralis\'{e}e $(x_{k}^{\ast })_{k\in I}$
dans la boule unit\'{e} de $(A_{0}^{\ast },A_{1}^{\ast })_{\beta }$ telle
que $x_{k}^{\ast }\rightarrow P^{\ast }x^{\ast },$ $\sigma ((Z^{\beta
})^{\ast },Z^{\beta })$ (on peut choisir la suite $(x_{k}^{\ast })_{k\in I}$
dans $A_{0}^{\ast }\cap A_{1}^{\ast }=B_{0}\cap B_{1}).$

Il est clair que pour tout $k\in I$ 
\begin{equation*}
-i(R^{\beta }\circ g(\beta +i\tau +i/n)-R^{\beta }\circ g(\beta +i\tau
),x_{k}^{\ast })=\left[ \dint\limits_{\tau }^{\tau +1/n}(R^{\beta }\circ
\phi _{\beta }(t),x_{k}^{\ast })dt\right] ;
\end{equation*}%
ceci implique que

\begin{eqnarray}
-i(R^{\beta }\circ (g(\beta +i\tau +i/n)-R^{\beta }\circ g(\beta +i\tau
),x_{{}}^{\ast }) &=&  \label{q} \\
-i(P\left[ R^{\beta }\circ (g(\beta +i\tau +i/n)\right] -P\left[ R^{\beta
}\circ g(\beta +i\tau )\right] ,x_{{}}^{\ast }) &=&  \notag \\
-i(R^{\beta }\circ (g(\beta +i\tau +i/n)-R^{\beta }\circ g(\beta +i\tau
),P^{\ast }x_{{}}^{\ast }) &=&  \notag \\
lim_{k}-i((R^{\beta }\circ g(\beta +i\tau +i/n)-R^{\beta }\circ g(\beta
+i\tau ),x_{k}^{\ast }) &=&.  \notag \\
&&lim_{k}\left[ \dint\limits_{\tau }^{\tau +1/n}(R^{\beta }\circ \phi
_{\beta }(t),x_{k}^{\ast })dt\right] .  \notag
\end{eqnarray}

D'apr\`{e}s l'\'{e}tape 1, $B_{F^{\ast }}$ est un compact de Rosenthal pour
la topolgie pr\'{e}faible, en appliquant le r\'{e}sultat de \cite{ROS}, on
voit que

\begin{eqnarray}
lim_{k}\left[ \dint\limits_{\tau }^{\tau +1/n}(R^{\beta }\circ \phi _{\beta
}(t),x_{k}^{\ast })dt\right] &=&  \label{r} \\
\left[ \dint\limits_{\tau }^{\tau +1/n}(\left[ R^{\beta }\circ \phi _{\beta
}(t)\right] ,P^{\ast }x_{{}}^{\ast })dt\right] . &=&  \notag \\
&&\dint\limits_{\tau }^{\tau +1/n}(P\left[ R^{\beta }\circ \phi _{\beta }(t)%
\right] ,x^{\ast })dt.  \notag
\end{eqnarray}

Donc, d'apr\`{e}s (\ref{q}) et (\ref{r}) on a

\begin{eqnarray*}
-i((R^{\beta }\circ g(\beta +i\tau +i/n)-R^{\beta }\circ g(\beta +i\tau
),x_{{}}^{\ast }) &=& \\
&&\left[ \dint\limits_{\tau }^{\tau +1/n}(P\left[ R^{\beta }\circ \phi
_{\beta }(t)\right] ,x_{{}}^{\ast })dt\right] .
\end{eqnarray*}

\emph{\'{e}tape 5:}$R^{\beta }\circ \phi _{\beta }=\phi _{1},$
presque-partout.

Nous appliquons maintenant le th\'{e}or\`{e}me de la differ\'{e}ntiabilit%
\'{e} de Labesgue \cite{DU} on obtient alors,%
\begin{equation}
n\left[ \dint\limits_{\tau }^{\tau +1/n}\phi _{1}(t)dt\right] \underset{%
n\rightarrow +\infty }{\rightarrow }\phi _{1}(\tau ),\text{ pour presque
tout }\tau \in \mathbb{R}.  \label{s}
\end{equation}

D'apr\`{e}s L'\'{e}tape 4, et la relation (\ref{p}), pour tout $n\in \mathbb{%
N}^{\ast }$ fix\'{e} et tout $x^{\ast }\in B_{0}^{\ast }\cap B_{1}^{\ast }$
on a

\begin{eqnarray*}
-in(R^{\beta }\circ g(\beta +i\tau +i/n)-R^{\beta }\circ g(\beta +i\tau
),x^{\ast }) &=& \\
n\left[ \dint\limits_{\tau }^{\tau +1/n}(P\left[ R^{\beta }\circ \phi
_{\beta }(t)\right] ,x^{\ast })dt\right] &=& \\
n\left[ \dint\limits_{\tau }^{\tau +1/n}(\phi _{1}(t),x^{\ast })dt\right] &=&
\\
n(\left[ \dint\limits_{\tau }^{\tau +1/n}\phi _{1}(t)dt\right] ,x^{\ast }),%
\text{ } &&
\end{eqnarray*}

c'est-\`{a}-dire%
\begin{eqnarray}
-in(R^{\beta }\circ g(\beta +i\tau +i/n)-R^{\beta }\circ g(\beta +i\tau
),x^{\ast }) &=&  \label{u} \\
&&n(\left[ \dint\limits_{\tau }^{\tau +1/n}\phi _{1}(t)dt\right] ,x^{\ast }).
\notag
\end{eqnarray}

Mais pour tout $\tau \in \mathbb{R}$ 
\begin{equation}
-in(R^{\beta }\circ g(\beta +i\tau +i/n)-R^{\beta }\circ g(\beta +i\tau
),x^{\ast })\underset{n\rightarrow +\infty }{\rightarrow }(R^{\beta }\circ
g^{\prime }(\beta +i\tau ),x^{\ast }),  \label{v}
\end{equation}%
car la fonction $R^{\beta }\circ g$ est holomorphe \`{a} valeurs dans $%
B_{0}^{\ast }+B_{1}^{\ast }.$

et d'apr\`{e}s (\ref{s})$,$ pour presque tout $\tau \in \mathbb{R},$ on a%
\begin{equation}
n(\left[ \dint\limits_{\tau }^{\tau +1/n}\phi _{1}(t)dt\right] ,x^{\ast })%
\underset{n\rightarrow +\infty }{\rightarrow }(\phi _{1}(\tau ),x^{\ast }),%
\text{ }\forall x^{\ast }\in A_{0}^{\ast }\cap A_{1}^{\ast }.  \label{t}
\end{equation}

\bigskip

\bigskip Finalement d'ap\`{e}s les relations (\ref{u}), (\ref{v}) et (\ref{t}%
) nous d\'{e}duisons que pour presque tout $\tau \in \mathbb{R}$ 
\begin{equation*}
(R^{\beta }\circ g^{\prime }(\beta +i\tau ),x^{\ast })=(\phi _{1}(\tau
),x^{\ast }),\forall x^{\ast }\in A_{0}^{\ast }\cap A_{1}^{\ast }
\end{equation*}%
et donc $R^{\beta }\circ \phi _{\beta }=\phi _{1},$ presque-partout.$%
\blacksquare $

\begin{corollary}
\label{nd}Soit $(A_{0},A_{1})$ un couple d'interpolation; supposons
qu\textquotedblright il existe un $\beta \in \left] 0,1\right[ $ tel que $%
A^{\beta }=A_{\beta }$ et que $A_{\beta }$ est un espace faiblement de Lindel%
\"{o}f. Alors $A^{\theta }=A_{\theta },$ pour tout $\theta \in \left] 0,1%
\right[ .$
\end{corollary}

\emph{Preuve:}Soit $g\in \mathcal{G}(A_{0},A_{1})$ tel que $g^{\prime
}(\theta )=a.$ D'apr\`{e}s la proposition \ref{gf}, il suffit de montrer que
l'application $\tau \in \mathbb{R}\rightarrow $

$(g^{\prime }(\beta +i\tau ),a^{\ast })\in A_{\beta }$ est de premi\`{e}re
classe de Baire, pour tout $a^{\ast }\in A_{\beta }^{\ast }.$

Soient $a^{\ast }\in A_{\beta }^{\ast }=(A_{0}^{\ast },A_{1}^{\ast })^{\beta
}$ et $h\in \mathcal{G}(A_{0}^{\ast },A_{1}^{\ast })$ tels que $h^{\prime
}(\beta )=a^{\ast };$ pour tout $n\in \mathbb{N}^{\ast }$ consid\'{e}rons $%
F_{n}(z)=-in\left[ h(z+i/n)-h(z)\right] ,$ $z\in S.$ Il est clair que la
suite $(F_{n}(\beta ))_{n\geq 0}$ est une suite born\'{e}e dans $%
(A_{0}^{\ast },A_{1}^{\ast })_{\beta }.$ D'autre part l'application $\tau
\in \mathbb{R}\rightarrow (\phi _{\beta }(\tau ),b^{\ast })$ est continue
pour tout $b^{\ast }\in A_{0}^{\ast }\cap A_{1}^{\ast },$ comme $A_{0}^{\ast
}\cap A_{1}^{\ast }$ est dense dans $(A_{0}^{\ast },A_{1}^{\ast })_{\beta }$%
, alors pour tout $n\in \mathbb{N}^{\ast },$ l'application $\tau \in \mathbb{%
R}\rightarrow (\phi _{\beta }(\tau ),F_{n}(\beta ))$ est continue.

Remarquons que pour tout $b^{{}}\in A_{0}^{{}}\cap A_{1}^{{}};$ $%
(b,F_{n}(\beta ))\underset{n\rightarrow +\infty }{\rightarrow }(b,h^{\prime
}(\beta )),$ donc $(\phi (\tau ),F_{n}(\beta ))\underset{n\rightarrow
+\infty }{\rightarrow }(\phi _{\beta }(\tau ),a^{\ast })$ (car $A_{0}\cap
A_{1}$ est dense dans $A_{\beta }=A^{\beta }).$ On en d\'{e}duit que
l'application$:\tau \in \mathbb{R}\rightarrow (g^{\prime }(\beta +i\tau
),a^{\ast })$ est de premi\`{e}re classe de Baire.$\blacksquare $

\begin{proposition}
\label{xc}Soit $(A_{0},A_{1})$ un couple d'interpolation; supposons qu'il
existe un $\beta \in \left] 0,1\right[ $ tel que
\end{proposition}

1)$(A_{0},A_{1})_{\beta }$ est faiblement s\'{e}quentiellement complet.

2$)(A_{0}^{\ast },A_{1}^{\ast })^{\beta }=(A_{0}^{\ast },A_{1}^{\ast
})_{\beta }.$

Alors $A^{\theta }=A_{\theta },$ pour tout $\theta \in \left] 0,1\right[ .$

Preuve:Soit $g\in \mathcal{G}(A_{0},A_{1})$; il suffit de montrer que
l'application $\phi _{\beta }:\tau \in \mathbb{R}\rightarrow g^{\prime
}(\beta +i\tau )$ est fortement mesurable.

Pour tout $n\in \mathbb{N}^{\ast },$ notons $F_{n}(\beta +i\tau )=-in\left[
g(\beta +i\tau +i/n)-g(\beta +i\tau )\right] .$ Fixons $\tau \in \mathbb{R};$
remarquons que la suite $(F_{n}(\beta +i\tau ))_{n\geq 0}$ est born\'{e}e
dans $A_{\beta }$ et que $\left[ (F_{n}(\beta +i\tau ),a^{\ast })\right]
_{n\geq 1}$ est de Cauchy, pour tout $a^{\ast }\in A_{0}^{\ast }\cap
A_{1}^{\ast },$ comme $A_{0}^{\ast }\cap A_{1}^{\ast }$ est dense dans $%
(A_{0}^{\ast },A_{1}^{\ast })_{\beta }$, d'apr\`{e}s la condition $(2),$ la
suite $(F_{n}(\beta +\iota \tau ))_{n\geq 1}$ est une faiblement de Cauchy
dans $A_{\beta }.$

D'apr\`{e}s la condition $(1),$ $(F_{n}(\beta +i\tau ))_{n\geq 1}$ est
converge vers $\phi _{\beta }(\tau )$ dans $A_{\beta }$ quand $n\rightarrow
+\infty $ (car $F_{n}(\beta +i\tau )\underset{n\rightarrow +\infty }{%
\rightarrow }\phi _{\beta }(\tau )$ dans $A_{0}+A_{1}),$ ce qui implique que 
$\phi _{\beta }$ est fortmement mesurable.$\blacksquare $

\begin{definition}
\label{cxw}Soit $(A_{0},A_{1})$ un couple d'interplation; on dit que le
couple d'interpolation $(A_{0},A_{1})$ a la propri\'{e}t\'{e} de
l'approximation \`{a} gauche (resp. \`{a} droite), s'il exist une suite g%
\'{e}n\'{e}ralis\'{e}e $(T_{i})_{i\in I}$ d'op\'{e}rateurs de $A_{0}+A_{1}$ 
\`{a} valeurs dans $A_{0}+A_{1}$ v\'{e}rifiant les propri\'{e}t\'{e}s
suivantes:
\end{definition}

I)Pour tout $i\in I,$ $T_{i}\in \mathcal{L}(A_{j}),$ $j=0,1$

II)Pour tout $i\in I,$ il existe $j\in \left\{ 0,1\right\} $ tel que $%
T_{i_{\mid _{A_{j}}}}$ est du rang fini.

III)$T_{i}\underset{}{\rightarrow }I_{A_{0}}$ dans $A_{0}$ uniform\'{e}ment
sur tout compact de $A_{0}$ (resp.$T_{i}\underset{}{\rightarrow }I_{A_{1}}$
dans $A_{1}$ uniform\'{e}ment sur tout compact d $A_{1}$).

IV) Pour tout compact $L$ de $A_{1},$ \ $C_{L}=\sup \left\{ \left\Vert
T_{i}x\right\Vert _{A_{1}};\text{ }i\in I\text{ et }x\in L\right\} <+\infty $
(resp. pour tout compact $L$ de $A_{0},$ \ $C_{L}=\sup \left\{ \left\Vert
T_{i}x\right\Vert _{A_{0}};\text{ }i\in I\text{ et }x\in L\right\} <+\infty
. $

\begin{definition}
\label{ik}Le couple d'interpolation $(A_{0},A_{1})$ est dit a la propri\'{e}t%
\'{e} de l'approximation, s'il a la propri\'{e}t\'{e} de l'approximation 
\`{a} gauche ou \`{a} droite.
\end{definition}

\begin{description}
\item 
\begin{proposition}
\label{bfm}Soient $(A_{0},A_{1})$ un couple d'interpolation qui a la propri%
\'{e}t\'{e} de l'approximation et $\theta \in \left] 0,1\right[ .$ Alors $%
A_{\theta }$ a la propri\'{e}t\'{e} de l'approximation \cite[chap.III-3,p.238%
]{DU}.
\end{proposition}
\end{description}

D\'{e}monstration.

Supposons que $(A_{0},A_{1})$ a la propri\'{e}t\'{e} de l'approximation \`{a}
gauche.

D'apr\`{e}s \cite[th.4.1.2]{BL} et la propri\'{e}t\'{e} (I)$,$ $T_{i}\in 
\mathcal{L}(A_{\theta })$ et d'apr\`{e}s la propri\'{e}t\'{e} (II), le rang
de $T_{i}$ est fini (comme op\'{e}rateur de $A_{\theta }$ \`{a} valeurs dans 
$A_{\theta }).$

Soient $a\in A_{\theta }$ et $F\in \mathcal{F(}A_{0},A_{1})$ tels que $%
F(\theta )=a;$ remarquons que $T_{i}\circ F\in \mathcal{F}(\overline{A})$ .
D'autre part le lemme 4.3.2 de \cite{BL} nous dit que 
\begin{eqnarray*}
\left\Vert T_{i}a-a\right\Vert _{A_{\theta }} &\leq &\left[ \dint\limits_{%
\mathbb{R}}\left\Vert T_{i}F(i\tau )-F(i\tau )\right\Vert
_{A_{0}}Q_{0}(\theta ,\tau )d\tau \right] ^{1-\theta }\times \\
&&\left[ \dint\limits_{\mathbb{R}}\left\Vert T_{i}F(1+i\tau )-F(1+i\tau
)\right\Vert _{A_{1}}Q_{1}(\theta ,\tau )d\tau \right] ^{\theta }.
\end{eqnarray*}

\bigskip Comme les ensembles $L=\left\{ F(i\tau );\text{ }\tau \in \mathbb{R}%
\right\} ,$ $L^{\prime }=\left\{ F(1+i\tau );\text{ }\tau \in \mathbb{R}%
\right\} $ sont relativement compacts de $A_{0}$ et $A_{1}$ respectivement,
d'apr\`{e}s les propri\'{e}t\'{e}s (III) et IV on a 
\begin{equation*}
\left[ \dint\limits_{\mathbb{R}}\left\Vert T_{i}F(i\tau )-F(i\tau
)\right\Vert _{A_{0}}Q_{0}(\theta ,\tau )d\tau \right] ^{1-\theta
}\rightarrow 0
\end{equation*}

\bigskip et%
\begin{equation*}
\left[ \dint\limits_{\mathbb{R}}\left\Vert T_{i}F(i\tau )-F(i\tau
)\right\Vert _{A_{1}}Q_{1}(\theta ,\tau )d\tau \right] ^{\theta }\leq \left[
C_{L^{\prime }}+\left\Vert F\right\Vert _{\mathcal{F}(\overline{A})}\right]
^{\theta }.
\end{equation*}

Il en resulte que $\left\Vert T_{i}a-a\right\Vert _{A_{\theta }}\rightarrow
0.\blacksquare $

\QTP{Body Math}
$\bigskip $

\QTP{Body Math}
$\bigskip $

\QTP{Body Math}
$\bigskip $

\end{document}